\documentclass[11pt]{article}
\usepackage{CJK}
\usepackage{amsmath,amsfonts,mathrsfs,amssymb}
\usepackage{indentfirst}
\numberwithin{equation}{section}
\usepackage{color}

\numberwithin{equation}{section}

\newcommand{\HT}{\CJKfamily{hei}}

\def\lb{\label}

\def\be{\begin{equation}}
\def\ee{\end{equation}}
\def\bea{\begin{eqnarray}}
\def\eea{\end{eqnarray}}
\def\bes{\begin{eqnarray*}}
\def\ees{\end{eqnarray*}}

\def\y{\begin{eqnarray*}}
\def\ey{\end{eqnarray*}}

\begin{document}
\title{\HT { Concentration-compactness principle of singular Trudinger-Moser inequality involving $N$-Finsler--Laplacian operator }}
\author{\small {Yanjun Liu }\\
\small $^{a}$School of Mathematical Sciences, Nankai University, Tianjin 300071, P. R. China}

\date{}

\maketitle \footnote[0] {Email: liuyj@mail.nankai.edu.cn(Y. Liu).}

\noindent{\small
{\bf Abstract:} In this paper, suppose $F: \mathbb{R}^{N} \rightarrow [0, +\infty)$ be a  convex function of class $C^{2}(\mathbb{R}^{N} \backslash \{0\})$ which is even and positively homogeneous of degree 1. We establish the Lions type concentration-compactness principle of singular Trudinger-Moser Inequalities involving $N$-Finsler--Laplacian operator. Let $\Omega\subset \mathbb{R}^{N}(N\geq 2)$ be a smooth bounded domain. $\{u_n\}\subset W_0^{1, N}(\Omega)$ be a sequence such that anisotropic Dirichlet norm$\int_{\Omega}F^N (\nabla u_n)dx=1$, $u_n \rightharpoonup u \not \equiv 0$ weakly in $W_0^{1, N}(\Omega)$. Then for any
 $0 < p < p_N(u):=(1-\int_{\Omega}F^N (\nabla u)dx)^{-\frac{1}{N-1}},$
we have
  $$ \int_{\Omega}\frac{e^{\lambda_{N}(1-\frac{\beta}{N})p |u_n|^{\frac{N}{N-1}}}}{F^{o}(x)^{\beta}}dx<+\infty,  $$
where $0\leq\beta <N$, $\lambda_{N}=N^{\frac{N}{N-1}} \kappa_{N}^{\frac{1}{N-1}}$ and $\kappa_{N}$ is the volume of a unit Wulff ball. This conclusion fails if $p \geq p_N(u)$.
 Furthermore, we also obtain the corresponding
concentration-compactness principle  in the entire Euclidean space $\mathbb{R}^{N}$.

}

\section{Introduction and main results}
 This paper is concerned with concentration-compactness-principle of singular Trudinger-Moser inequality involving N-Finsler-Laplacian operator. In order to give our motivation, let's recall some known results. Suppose $\Omega \subset \mathbb{R}^{N}(N \geq 2)$ be a bounded smooth domain.
 $W_{0}^{1,N}(\Omega)\hookrightarrow L^q(\Omega)$ for $1 \leq q < \infty$, but the embedding  $W_{0}^{1,N}(\Omega)\not \hookrightarrow L^\infty(\Omega)$, one can see the counterexample by taking $u(x) = (- \ln| \ln |x||)_{+}$ as $\Omega$ is the unit ball. It was proposed independently by Yudovich \cite{Yudovich}, Pohozaev \cite{Pohozaev},  Peetre \cite{Peetre} and Trudinger \cite{N.S. Trudinger} that $W^{1, N}_{0}(\Omega)$ is embedded in the Orlicz space $L_{\varphi_{\alpha}}(\Omega)$ determined by the Young function $\varphi_{\alpha}(t)=e^{\alpha|t|^{\frac{N}{N-1}}}-1$ for some positive number $\alpha$. Moser \cite{J. Moser} sharpened the results of Trudinger \cite{N.S. Trudinger} and established the following  inequality
\begin{equation}\label{1.1}
 \sup_{u\in W_0^{1, N}(\Omega), \|\nabla u\|_{N}\leq1}\int_{\Omega}e^{\alpha|u|^{\frac{N}{N-1}}}dx<+\infty,  \forall \alpha \leq  \alpha_{N},
\end{equation}
where $\alpha_{N}:=N\omega_{N-1}^{1/(N-1)}$, $\omega_{N-1}$ is the surface measure of the
unit sphere in $\mathbb{R}^{N}$.
Moreover, the supremum in ($\ref{1.1}$) is $+\infty$ if $\alpha > \alpha_N$. Inequality ($\ref{1.1}$) is now referred as Trudinger-Moser inequality and plays an important role in geometric analysis and partial differential equations (see \cite{ddR}).   Using a rearrangement argument and a change of variables, Adimurthi-Sandeep \cite{AS} generalized the Trudinger-Moser inequality to a singular version as follows:
\be   \sup_{u\in W_0^{1, N}(\Omega), \int_{\Omega}|\nabla u|^{N}dx\leq 1}\int_{\Omega}\frac{e^{\alpha u^{\frac{N}{N-1}}}}{|x|^{\beta}}dx<+\infty,   \lb{1.2}\ee
where $0\leq \beta< N$, $0<\alpha \leq \alpha_{N}(1-\frac{\eta}{N})$, $\alpha_{N}=N\omega_{N-1}^{\frac{1}{N-1}}$, $\omega_{N-1}$ is the surface measure of the unit sphere in $\mathbb R^N$. Moreover, this inequality is sharp, {\it i.e.},  when $\alpha > \alpha_{N}(1-\frac{\eta}{N})$, the supremum is infinity. Trudinger-Moser inequalities for unbounded domains were proposed by D. M. Cao \cite{D.M. Cao} in dimension two and J. M. do \'O \cite{J. M. do}, Adachi-Tanaka \cite{AT} in high dimension. Ruf \cite{Ruf} (for the case $N=2$), Li and Ruf  \cite{LB}  (for the general case $N \geq 2$) obtained the Trudinger-Moser inequality in the critical case by replacing the Dirichlet norm with the standard Sobolev norm in $W^{1, N}(\mathbb{R}^{N})$. Obviously, if $\beta = 0$, then \eqref{1.2} reduces to the famous Trudinger-Moser inequality. Subsequently, the inequality \eqref{1.2} was extended to the entire Euclidean space $\mathbb{R}^{N}$ by Adimurthi-Yang \cite{AY}.

An important result is concentration-compactness principle  with Trudinger-Moser inequality due to P. L Lions \cite{Lions}. More precisely, let $\{u_n\}\subset W_0^{1, N}(\Omega)$ be a sequence such that $\| \nabla u_n \|_N=1$, $u_n \rightharpoonup u \not \equiv 0$ weakly in $W_0^{1, N}(\Omega)$. Then for any $0 < p < p_N:=(1-\int_\Omega |\nabla u|^N)^{-\frac{1}{N-1}}$,
it holds
\begin{equation}\label{4}
   \int_{\Omega}e^{\alpha_N p |u_n|^{\frac{N}{N-1}}}dx<+\infty .
\end{equation}
 Roughly speaking, the concentration-compactness principle  tells us that, if a sequence $\{u_n\}\subset W_0^{1, N}(\Omega)$ converges  weakly to some function $u \in W_0^{1, N}(\Omega)$, and does not concentrate at one point in $\Omega$, then
an inequality like ($\ref{4}$) holds along the sequence $\{u_n\}$, with a constant larger than $N\omega_{N-1}^{1/(N-1)}$,
depending on $\| \nabla u\|_{N}$. In \cite{Lions}, the author only proved the case of $0<p<p^*_N:==(1-\int_{\Omega{^*}} |\nabla u^*|^N)^{-\frac{1}{N-1}}$, we know $p^*_N\leq p_N$ by Poly\'a-Szeg\"o inequality $\int_{\Omega{^*}} |\nabla u^*|^Ndx\leq\int_\Omega |\nabla u|^Ndx$, here $u^*$ is the radially decreasing symmetry of $u$. We should pay  attention to the recent work in \cite{CCH} by \v{C}ern\'{y}  et al. The authors present a new proof of this relevant principle for $0<p<p_N$. Moreover, this approach allows one to treat functions
with unrestricted boundary values in bounded domains. Concentration-compactness principle is a powerful tool in proving existence of extremal functions and existence of solutions to boundary value problems. It has been extended to the singular version in \cite{GS}.  Their results can be stated as follows: let $\{u_n\}\subset W_0^{1, N}(\Omega)$ be a sequence such that $\| \nabla u_n \|_N=1$, $u_n \rightharpoonup u \not \equiv 0$ weakly in $W_0^{1, N}(\Omega)$, $\nabla u_n \rightharpoonup \nabla u$ a.e. in $\Omega$. Then for any $0\leq\beta <N$ and $0 < p < (1-\| \nabla u\|_N^N)^{-\frac{1}{N-1}}$,
there holds
\begin{equation}\label{LL}
   \int_{\Omega}\frac{e^{\alpha_N(1-\frac{\beta}{N}) p |u_n|^{\frac{N}{N-1}}}}{|x|^\beta}dx<+\infty .
\end{equation}
More concentration-compactness principle on unbounded domain and the Heisenberg group, we refer the reader to \cite{do,LiLu,Zhang}.

Another interesting research is that Trudinger-Moser inequality has been generalized to the case of anisotropic norm.  In this paper, denote that $F \in C^{ 2}(\mathbb{R}^{N} \backslash {0})$ is a positive, convex and homogeneous function, $F_{\xi_{i}} = \frac{\partial F}{\partial\xi_{i}}$ and its polar $F^{o}(x)$ represents a Finsler metric on $\mathbb{R}^{N}$. We will replace the isotropic Dirichlet norm $\|u\|_{W_0^{1, N}(\Omega)}=(\int_\Omega |\nabla u|^Ndx)^{\frac{1}{N}}$ by the anisotropic Dirichlet norm $(\int_\Omega F^N(\nabla u)dx)^{\frac{1}{N}}$ in $W_0^{1, N}(\Omega)$. In \cite{WX}, Wang and Xia proved the following result:\\
\textbf{Theorem A.} Suppose $\Omega\subset \mathbb{R}^{N}(N\geq 2)$ be a smooth bounded domain. Let $u\in W_0^{1, N}(\Omega)$ and $(\int_{\Omega}F^{N}(\nabla u)dx)\leq 1$. Then there exists a constant $C(N)$, such that
   \be \int_{\Omega}e^{\lambda u^{\frac{N}{N-1}}}dx  \leq C(N) |\Omega|, \lb{(1.4)}\ee
where $0<\lambda \leq \lambda_{N}=N^{\frac{N}{N-1}} \kappa_{N}^{\frac{1}{N-1}}$ and $\kappa_{N}=|\{x \in \mathbb{R}^{N} : F^{o}(x)\leq 1\}|$. $\lambda_{N}$ is sharp in the sense that if $\lambda> \lambda_{N}$ then there exists a sequence $(u_{n})$ such that $\int_{\Omega}e^{\lambda u^{\frac{N}{N-1}}}dx $ diverges. In \cite{ZZ}, the authors obtained the existence of extremal functions for the sharp geometric inequality \eqref{(1.4)}.

For the minimization problem of $\int_{\Omega}F^{N}(\nabla u)dx$, we know that its Euler equation contains an operator of the form
   $$Q_{N}u:=\sum_{i=1}^{N}\frac{\partial}{\partial x_{i}}(F^{N-1}(\nabla u)F_{\xi_{i}}(\nabla u)),$$
which is called $N$-Finsler-Laplacian operator. When $N=2$ and $F(\xi) = |\xi|$, $Q_{2}$ is just the ordinary Laplacian. The operator $Q_N$ is closely related to a smooth, convex hypersurface in $\mathbb{R}^{N}$. It has been studied in some literatures, see \cite{AFTL,BF,FK} and the references therein.  We denote $\kappa_{N}=|\{x \in \mathbb{R}^{N} : F^{o}(x)\leq 1\}|$ is the volume of a unit Wulff ball. Recently, by using a convex symmetrization approach proposed in  \cite{AFTL},  which is the extension of Schwarz
symmetrization in \cite{T}, X. Zhu \cite{Zhu} derived the following results.\\
\textbf{Theorem 1.1.} (see \cite{Zhu}) Let $\Omega\subset \mathbb{R}^{N}(N\geq 2)$ be a smooth bounded domain.  Then
  \be   \sup_{u\in W^{1, N}_{0}(\Omega), \int_{\Omega}F^{N}(\nabla u)dx\leq 1}\int_{\Omega}\frac{e^{\lambda |u|^{\frac{N}{N-1}}}}{F^{o}(x)^{\beta}}dx<+\infty,   \lb{(1.5)}\ee
and
   \be   \sup_{u\in W^{1, N}(\mathbb{R}^{N}), \int_{\mathbb{R}^{N}}(F^{N}(\nabla u)+\tau |u|^N)dx\leq 1}\int_{\mathbb{R}^{N}}\frac{\Phi(\lambda |u_{n}|^{\frac{N}{N-1}})}{F^{o}(x)^{\beta}}dx<\infty,      \lb{(1.6)}\ee
where $0\leq \beta< N$, $\tau>0$, $\Phi(s):=e^{s}-\sum_{k=0}^{N-2}\frac{s^{k}}{k!}$, $0<\lambda \leq \lambda_{N}(1-\frac{\beta}{N})$, $\lambda_{N}=N^{\frac{N}{N-1}} \kappa_{N}^{\frac{1}{N-1}}$ and $\kappa_{N}$ is the volume of a unit Wulff ball.  Moreover, the above inequalities are sharp, {\it i.e.},  when $\lambda > \lambda_{N}(1-\frac{\beta}{N})$, the supremum is infinity.

In this paper, we will establish the Lions type concentration-compactness principle of singular Trudinger-Moser Inequalities under the anisotropic norm.   \\
\textbf{Theorem 1.2 } Let $\Omega\subset \mathbb{R}^{N}(N\geq 2)$ be a smooth bounded domain. $\{u_n\}\subset W_0^{1, N}(\Omega)$ be a sequence such that $\int_{\Omega}F^N (\nabla u_n)dx=1$, $u_n \rightharpoonup u \not \equiv 0$ weakly in $W_0^{1, N}(\Omega)$. Then for any
     $$0 < p < p_N(u):=(1-\int_{\Omega}F^N (\nabla u)dx)^{-\frac{1}{N-1}},$$
we have
  \be \int_{\Omega}\frac{e^{\lambda_{N}(1-\frac{\beta}{N})p |u_n|^{\frac{N}{N-1}}}}{F^{o}(x)^{\beta}}dx<+\infty  \lb{(2.2)}\ee
where $0\leq \beta <N$, $\lambda_{N}=N^{\frac{N}{N-1}} \kappa_{N}^{\frac{1}{N-1}}$ and $\kappa_{N}$ is the volume of a unit Wulff ball. Moreover, this conclusion fails if $p \geq p_N(u)$.\\
\textbf{Theorem 1.3 } Suppose $\{u_n\}\subset W^{1, N}(\mathbb{R}^{N})$ be a sequence such that $\int_{\mathbb{R}^{N}}(F^N (\nabla u_n)+|u_n|^N)dx=1$, $u_n \rightharpoonup u \not \equiv 0$ weakly in $W^{1, N}(\mathbb{R}^{N})$. Then for any
     $$0 < p < \bar{p}_N(u):=(1-\int_{\mathbb{R}^{N}}(F^N (\nabla u)+|u|^N)dx)^{-\frac{1}{N-1}},$$
we have
  \be \int_{\mathbb{R}^{N}}\frac{\Phi(\lambda_{N}(1-\frac{\beta}{N})p |u_n|^{\frac{N}{N-1}})}{F^{o}(x)^{\beta}} dx<+\infty  \lb{(2.2)}\ee
where$0\leq \beta <N$, $\Phi(s):=e^{s}-\sum_{k=0}^{N-2}\frac{s^{k}}{k!}$, $\lambda_{N}=N^{\frac{N}{N-1}} \kappa_{N}^{\frac{1}{N-1}}$ and $\kappa_{N}$ is the volume of a unit Wulff ball. Moreover, this conclusion fails if $p \geq \bar{p}_N(u)$.

This paper is organized as follows: In Section 2, we give some
preliminaries. In Section 3, we establish the Lions type concentration-compactness principle of
singular Trudinger-Moser Inequality under the anisotropic Direchlet norm. In Section 4, we  obtain the corresponding
concentration-compactness principle  in the entire Euclidean space $\mathbb{R}^{N}$.
\section{preliminaries }
In this section, we will give some preliminaries for our use later.

Let $F: \mathbb{R}^{N} \rightarrow [0, +\infty)$ be a  convex function of class $C^{2}(\mathbb{R}^{N} \backslash \{0\})$ which is even and positively homogeneous of degree 1, so that
\be    F(t \xi)=|t|F(\xi)~~~~~\operatorname{for~~any}~~t\in \mathbb{R},~~ \xi \in \mathbb{R}^{N}. \lb{2.1}\ee
We also  assume that $  F(\xi)>0$ for any $\xi \neq 0$  and  $Hess(F^2)$ is positive definite in $\mathbb{R}^{N} \backslash \{0\}$. A typical example is $F(\xi)=(\sum_{i}|\xi_i|^q)^{\frac{1}{q}}$ for $q\in [1, \infty)$.

Let $F^{o}$ be the support function of $K:=\{x\in \mathbb{R}^{N}: F(x) \leq 1 \}$, which is defined by
      $$F^{o}(x):= \sup_{\xi \in K}\langle x, \xi \rangle, $$
so $F^{o}: \mathbb{R}^{N} \rightarrow [0, +\infty)$ is also a convex, homogeneous function of class  $C^{2}(\mathbb{R}^{N} \backslash \{0\})$. From \cite{AFTL}, $F^{o}$ is dual to $F$ in the sense that
       $$F^{o}(x)= \sup_{\xi \neq 0}\frac{\langle x, \xi \rangle}{F(\xi)},~~~~F(x)= \sup_{\xi \neq 0}\frac{\langle x, \xi \rangle}{F^o(\xi)}.$$

Consider the map $\phi: S^{N-1} \rightarrow \mathbb{R}^{N}$, $\phi(\xi)=F_\xi(\xi)$. Its image $\phi(S^{N-1})$
is smooth, convex hypersurface in $\mathbb{R}^{N}$, which is called the Wulff shape (or equilibrium crystal shape) of $F$. Then $\phi(S^{N-1})= \{x\in \mathbb{R}^{N}|F^o(x) = 1\}$(see \cite{WX1}, Proposition 2.1).

We also give some simple properties of the function $F$, which follows directly from the assumption on $F$, also see \cite{WX,FK}.\\
\textbf{Lemma 2.1.} There hold\\
 (i) $|F(x)-F(y)|\leq F(x+y) \leq F(x)+F(y)$;\\
 (ii) $\frac{1}{C}\leq |\nabla F(x)|\leq C$ and  $\frac{1}{C}\leq |\nabla F^{o}(x)|\leq C$ for some $C>0$ and any $x\neq 0$; \\
 (iii) $\langle x, \nabla F(x)\rangle= F(x), \langle x, \nabla F^{o}(x)\rangle = F^{o}(x)$ for  any $x\neq 0$.\\
\textbf{Remark 2.2.} Since $Hess(F^2)$ is positive definite in $\mathbb{R}^{N} \backslash \{0\}$. Then by
Xie and Gong \cite{XG}, $Hess(F^ N)$ is also positive definite in $\mathbb{R}^{N} \backslash \{0\}$. Moreover, for a bounded smooth domain $\Omega \subset \mathbb{R}^{N}(N \geq 2)$, we know that $Q_2$ is a uniformly elliptic operator in any compact subsets of $\Omega \backslash \{x|\nabla u(x)=0\}$, see \cite{WX1}.

We will use the convex symmetrization which is defined in \cite{AFTL}. The convex symmetrization generalizes the Schwarz symmetrization(see \cite{T}). Let us consider a measured function $u$ on $\Omega \subset \mathbb{R}^{N}$, one dimensional decreasing rearrangement of $u$ is
   \be u^{\sharp}(t)=\sup\{s\geq 0: |\{x \in \Omega:   |u(x)|> s \}|>t \}~~~\operatorname{for}~~t\in  \mathbb{R}. \lb{2.2}\ee
The convex symmetrization of $u$ with respect to $F$ is defined as
   \be u^{\star}(x)=u^{\sharp}(\kappa_N F^{o}(x)^{N}) ~~~\operatorname{for}~~x\in \Omega^{\star}.  \lb{2.3}\ee
Here $\kappa_N F^{o}(x)^{N}$ is just the Lebesgue measure of a homothetic Wulff ball with radius $F^{o}(x)$ and $\Omega^{\star}$ is the homothetic Wulff ball centered at the origin having the same measure as $\Omega$. In \cite{AFTL}, the authors proved a P\'olya-Szeg\"o principle and a comparison result for solutions of the
Dirichlet problem for elliptic equations for the convex symmetrization, which generalizes the classical
results for Schwarz symmetrization due to Talenti \cite{T}.\\
\textbf{Lemma 2.3.} (see \cite{AFTL}) If $u\in W_0^{1, p}(\Omega)$ for $p\geq 1$. Then $u^{\star}\in  W_0^{1, p}(\Omega^{\star})$ and
   $$\int_{\Omega}F^p (\nabla u)dx\geq \int_{\Omega^{\star}}F^p (\nabla u^{\star})dx.$$

 Next, we denote $D_u(\mu)=\{x \in \Omega:  |u(x)|\geq\mu\}$. It is easily derived
   \be u^\star(x)=\sup\{\mu : F^o(x)\leq r, \kappa_N r^N=|D_u(\mu)|\}.     \lb{2.4}\ee
We claim: for any $p\in [0, N)$, it holds
   \be \bigg(\frac{1}{(F^o)^p}\bigg)^\star(x)\leq\frac{1}{F^o(x)^p}.     \lb{2.5}\ee
In fact, by \eqref{2.4}, we have
    \be \bigg(\frac{1}{(F^o)^p}\bigg)^\star(x)=\sup\{\mu : F^o(x)\leq r, \kappa_N r^N=|D_{(F^o)^{-p}}(\mu)|\}. \lb{2.6}\ee
According to our notation, we have
\begin{equation}
\left.
\begin{aligned}[b]
      D_{(F^o)^{-p}}(\mu)= &\{x \in \Omega:  (F^o(x))^{-p}\geq\mu\} \\
      = &\{x \in \Omega:  F^o(x)\leq \frac{1}{\mu^{1/p}}\} \\
       \subset   &\{x \in \mathbb{R}^{N}:  F^o(x)\leq \frac{1}{\mu^{1/p}}\}
\end{aligned}
\right.   \lb{(2.7)}
\end{equation}
Thus $|D_{(F^o)^{-p}}(\mu)|\leq \kappa_N \frac{1}{\mu^{N/p}}$, combing with \eqref{2.6}, we have $\kappa_Nr^N\leq \kappa_N \frac{1}{\mu^{N/p}}$, so $\mu\leq \frac{1}{r^p}$.
Therefore
\begin{equation}
\left.
\begin{aligned}[b]
 \bigg(\frac{1}{(F^o)^p}\bigg)^\star(x)=&\sup\{\mu : F^o(x)\leq r, \kappa_N r^N=|D_{(F^o)^{-p}}(\mu)|\}\\
                                       \leq &\inf\{\frac{1}{r^p} : F^o(x)\leq r \}=\frac{1}{(F^o(x))^p}
 \end{aligned}
\right.
\end{equation}
Our claim is proved.

 Now suppose $h$ and $\varphi$ be real-valued functions defined for $x\in \Omega$ with $h$ integrable over $\Omega$. Let $\varphi$ be measurable over $\Omega$ and satisfy the condition  $-\infty<\varphi_0\leq \varphi(x)\leq \varphi_1<\infty$, set  $D_\varphi(t)=\{x \in \Omega:  \varphi(x)\geq t\}$. Then Lemma 2.3 in \cite{Bandle} implies
     \be \int_\Omega h\varphi dx= \varphi_0\int_\Omega hdx+\int_{\varphi_0}^{\varphi_1}dt\int_{D_\varphi(t)}h dx.     \lb{2.8}\ee
\textbf{Lemma 2.4.} Assume that $f: [\varphi_0, \varphi_1]\rightarrow \mathbb{R}^{+}$ is a increasing function. Then we have
 \be \int_\Omega hf(\varphi) dx\leq \int_{\Omega^\star} h^\star f(\varphi^\star) dx.     \lb{2.}\ee
\textbf{Proof.} On one hand,
\begin{equation}
\left.
\begin{aligned}[b]
\int_\Omega hf(\varphi) dx=& f(\varphi_0)\int_\Omega hdx+\int_{f(\varphi_0)}^{f(\varphi_1)}dt\int_{\{x\in \Omega: f(\varphi)\geq t\}} h dx\\
                         =& f(\varphi_0)\int_\Omega hdx+\int_{f(\varphi_0)}^{f(\varphi_1)}dt\int_{\{x\in \Omega: \varphi\geq f^{-1}(t)\}} h dx.
\end{aligned}
\right.  \lb{(2.10)}
\end{equation}
On the other hand, since $\inf \varphi =\inf \varphi^\star$ and $\sup \varphi =\sup \varphi^\star$,
\begin{equation}
\left.
\begin{aligned}[b]
\int_{\Omega^\star} h^\star f(\varphi^\star) dx
=& f(\varphi_0)\int_{\Omega^\star} h^\star dx+\int_{f(\varphi_0)}^{f(\varphi_1)}dt\int_{\{x\in \Omega^\star: f(\varphi^\star)\geq t\}} h^\star dx\\
                         =& f(\varphi_0)\int_{\Omega^\star} h^\star dx+\int_{f(\varphi_0)}^{f(\varphi_1)}dt\int_{\{x\in \Omega^\star: \varphi^\star \geq f^{-1}(t)\}}  h^\star dx\\
                          =& f(\varphi_0)\int_{\Omega^\star} h^\star dx+\int_{f(\varphi_0)}^{f(\varphi_1)}dt\int_{\{x\in \Omega^\star: \varphi \geq f^{-1}(t)\}^\star}  h^\star dx.
\end{aligned}
\right.  \lb{(2.11)}
\end{equation}
Notice that $\int_\Omega hdx=\int_{\Omega^\star} h^\star dx$ and Lemma 2.2 in \cite{Bandle} implies
     \be \int_{\{x\in \Omega: \varphi\geq f^{-1}(t)\}} h dx\leq  \int_{\{x\in \Omega^\star: \varphi \geq f^{-1}(t)\}^\star}  h^\star dx.     \lb{2.12}\ee
The assertion now follows immediately.    $\hfill\Box$
\section{Lions type concentration-compactness principle in bounded domain}
In this section, we will prove Lions type concentration-compactness principle of
singular Trudinger-Moser Inequalities under the anisotropic Dirichlet norm, which can  be refered to \cite{LiLu} and Lemma 2.3 in \cite{ZZ}. This is the extension of Concentration-Compactness Principle due to P. L. Lions \cite{Lions}.\\
\textbf{Proof of Theorem 1.2.} From  the weak semicontinuity of the norm in $W^{1, N}_0(\Omega)$, we have
   $$\int_{\Omega}F^{N}(\nabla u)dx \leq \liminf_{n\rightarrow \infty}\int_{\Omega}F^{N}(\nabla u_n)dx=1.$$
Firstly, let $0<\int_{\Omega}F^N (\nabla u)dx<1 $, we give the proof by contradiction. Assume that there exists some $p_1 < p_N(u)$  and a subsequence of $\{u_n\}$(still denote $u_n$) such that
     \be \sup_n\int_{\Omega}\frac{e^{\lambda_{N}(1-\frac{\beta}{N})p_1 |u_n|^{\frac{N}{N-1}}}}{F^{o}(x)^{\beta}}dx=+\infty.  \lb{(2.3)}\ee
Set $\Omega_L^{n}=\{x\in \Omega : |u_n(x)|\geq L\}$, where $L$ is a positive constant. Let $v_n=|u_n|-L$, for any $\epsilon > 0$ and some positive constant $C$, by Young inequality $a^q b^{q'}\leq \epsilon a+\epsilon^{-\frac{q}{q'}}b, \frac{1}{q}+\frac{1}{q'}=1$  we have
\begin{equation}
\left.
\begin{aligned}[b]
  |u_n|^{\frac{N}{N-1}}\leq  & |v_n|^{\frac{N}{N-1}}+C|v_n|^{\frac{1}{N-1}}L +   L^{\frac{N}{N-1}}\\
                      =  & |v_n|^{\frac{N}{N-1}}+C(|v_n|^{\frac{N}{N-1}})^{\frac{1}{N}}(L^{\frac{N}{N-1}})^{\frac{N-1}{N}}+ L^{\frac{N}{N-1}}\\
                      \leq  & |v_n|^{\frac{N}{N-1}}+C\cdot\bigg(\frac{\epsilon}{C}|v_n|^{\frac{N}{N-1}}+(\frac{\epsilon}{C})^{-\frac{1}{N-1}}L^{\frac{N}{N-1}}\bigg)+ L^{\frac{N}{N-1}}\\
                       =: & (1+\epsilon)v_n^{\frac{N}{N-1}}+C_\epsilon L^{\frac{N}{N-1}}.
\end{aligned}
\right.  \lb{(2.4)}
\end{equation}
Since $0\leq \beta <N$, we have
 \begin{equation}
\left.
\begin{aligned}[b]
              \int_{\Omega}\frac{e^{\lambda_{N}(1-\frac{\beta}{N})p_1 |u_n|^{\frac{N}{N-1}}}}{F^{o}(x)^{\beta}}dx= & \int_{\Omega_L^{n}}\frac{e^{\lambda_{N}(1-\frac{\beta}{N})p_1 |u_n|^{\frac{N}{N-1}}}}{F^{o}(x)^{\beta}}dx+\int_{\Omega\backslash\Omega_L^{n}}\frac{e^{\lambda_{N}(1-\frac{\beta}{N})p_1 |u_n|^{\frac{N}{N-1}}}}{F^{o}(x)^{\beta}}dx\\
             \leq  & \int_{\Omega_L^{n}}\frac{e^{\lambda_{N}(1-\frac{\beta}{N})p_1 |u_n|^{\frac{N}{N-1}}}}{F^{o}(x)^{\beta}}dx+e^{\lambda_{N}(1-\frac{\beta}{N})p_1 L^{\frac{N}{N-1}}}\int_{\Omega}\frac{1}{F^{o}(x)^{\beta}}dx\\
              \leq  & \int_{\Omega_L^{n}}\frac{e^{\lambda_{N}(1-\frac{\beta}{N})p_1 |u_n|^{\frac{N}{N-1}}}}{F^{o}(x)^{\beta}}dx+C(L, N, \beta),
\end{aligned}
\right.  \lb{(2.5)}
\end{equation}
and then
\be \sup_n\int_{\Omega_L^{n}}\frac{e^{\lambda_{N}(1-\frac{\beta}{N})p_1 |u_n|^{\frac{N}{N-1}}}}{F^{o}(x)^{\beta}}dx=+\infty.  \lb{(2.6)}\ee
From \eqref{(2.4)}, we have
\begin{equation}
\left.
\begin{aligned}[b]
              \int_{\Omega_L^{n}}\frac{e^{\lambda_{N}(1-\frac{\beta}{N})p_1 |u_n|^{\frac{N}{N-1}}}}{F^{o}(x)^{\beta}}dx
             \leq  & e^{\lambda_{N}(1-\frac{\beta}{N})p_1 C _\epsilon L^{\frac{N}{N-1}}}
             \int_{\Omega_L^{n}}\frac{e^{(1+\epsilon)\lambda_{N}(1-\frac{\beta}{N})p_1 |v_n|^{\frac{N}{N-1}}}}{F^{o}(x)^{\beta}}dx.
\end{aligned}
\right.  \lb{(2.7)}
\end{equation}
Thus
\be \sup_n\int_{\Omega_L^{n}}\frac{e^{\lambda_{N}(1-\frac{\beta}{N})\overline{p}_1 |v_n|^{\frac{N}{N-1}}}}{F^{o}(x)^{\beta}}dx=\sup_n\int_{\Omega_L^{n}}\frac{e^{\lambda_{N}(1-\frac{\beta}{N})(\overline{p}_1^{\frac{N-1}{N}} |v_n|)^{\frac{N}{N-1}}}}{F^{o}(x)^{\beta}}dx=+\infty,  \lb{(2.8)}\ee
where $\overline{p}_1=(1+\epsilon)p_1<p_N(u)$. Now, we define
  $$T^L(u)=\min \{L, |u|\}sign(u) ~~~ and ~~~T_L(u)=u-T^L(u)$$
From the assumption  $0<\int_{\Omega}F^N (\nabla u)dx<1 $, we choose $L$  large enough such that
     \be \frac{1-\int_\Omega F^{N}(\nabla u)dx}{1-\int_\Omega F^{N}(\nabla T^L(u))dx} > \bigg(\frac{\overline{p}_1}{p_N(u)}\bigg)^{N-1} \lb{(2.9)}. \ee
Since $T^L(u_n)$ is  bounded in $W^{1, N}_0(\Omega)$, hence, up to a subsequence, $T^L(u_n)\rightharpoonup T^L(u)$ in $W^{1, N}_0(\Omega)$ and
$T^L(u_n)\rightarrow T^L(u)$ a.e. in $\Omega$. Combing \eqref{(2.8)} and \eqref{(1.5)}, up to a subsequence, we have
 $$  \limsup_{n\rightarrow \infty}\int_{\Omega_L^{n}}\overline{p}_1^{\frac{N-1}{N}}F^{N}(\nabla v_n)dx=\limsup_{n\rightarrow \infty}\int_{\Omega_L^{n}}F^{N}(\overline{p}_1^{\frac{N-1}{N}}\nabla v_n)dx\geq 1,  $$
which implies
 \be \int_{\Omega_L^{n}}F^{N}(\nabla v_n)dx= \int_{\Omega}F^{N}(\nabla T_L(u_n))dx\geq \bigg(\frac{1}{\overline{p}_1}\bigg)^{N-1}+o_n(1).  \lb{(2.10)} \ee
Thus,
\begin{equation*}
\left.
\begin{aligned}[b]
             &\bigg(\frac{1}{\overline{p}_1}\bigg)^{N-1}+\int_{\Omega}F^{N}(\nabla T^L(u_n))dx+o_n(1)\\
             \leq  & \int_{\Omega}F^{N}(\nabla T_L(u_n))dx+\int_{\Omega\backslash \Omega_L^{n}}F^{N}(\nabla u_n)dx\\
             =  & \int_{\Omega_L^{n}}F^{N}(\nabla u_n)dx+\int_{\Omega\backslash \Omega_L^{n}}F^{N}(\nabla u_n)dx=1.
\end{aligned}
\right.
\end{equation*}
The above inequality, the weak lower semicontinuity of norm, and \eqref{(2.9)} yield
\begin{equation*}
\left.
\begin{aligned}[b]
             \overline{p}_1\geq & \frac{1}{(1- \liminf_{n\rightarrow \infty}\int_{\Omega}F^{N}(\nabla T^L(u_n))dx )^{\frac{1}{N-1}}}\\
            \geq & \frac{1}{(1- \int_{\Omega}F^{N}(\nabla T^L(u))dx )^{\frac{1}{N-1}}}\\
             > & \frac{\overline{p}_1}{p_N(u)}   \frac{1}{(1- \int_{\Omega}F^{N}(\nabla u)dx )^{\frac{1}{N-1}}}=\overline{p}_1,
\end{aligned}
\right.
\end{equation*}
which is a contradiction.  Secondly, let $\int_{\Omega}F^N (\nabla u)dx=1 $, we can repeat the process of first case and get
  $$ \sup_n\int_{\Omega_L^{n}}\frac{e^{\lambda_{N}(1-\frac{\beta}{N})\overline{p}_1 |v_n|^{\frac{N}{N-1}}}}{F^{o}(x)^{\beta}}dx=+\infty, $$
where $\overline{p}_1=(1+\epsilon)p_1$. Then we have
    $$ \limsup_{n\rightarrow \infty}\int_{\Omega_L^{n}}F^{N}(\nabla v_n)dx =\limsup_{n\rightarrow \infty}\int_{\Omega}F^{N}(\nabla T_L(u_n))dx\geq \bigg(\frac{1}{\overline{p}_1}\bigg)^{N-1},$$
thus,
\begin{equation}
\left.
\begin{aligned}[b]
              \int_{\Omega}F^{N}(\nabla T^L(u))dx\leq &  \liminf_{n\rightarrow \infty} \int_{\Omega}F^{N}(\nabla T^L(u_n))dx\\
                                                  =&1- \limsup_{n\rightarrow \infty}\int_{\Omega}F^{N}(\nabla T_L(u_n))dx\\
                                                 \leq& 1-\bigg(\frac{1}{\overline{p}_1}\bigg)^{N-1}.
\end{aligned}
\right.  \lb{(2.12)}
\end{equation}
On the other hand, since $\int_{\Omega}F^N (\nabla u)dx=1$, we can choose $L>0$ in such a way that
\be \int_{\Omega}F^{N}(\nabla T^L(u))dx>1-\frac{1}{2} \bigg(\frac{1}{\overline{p}_1}\bigg)^{N-1}. \lb{(2.13)}\ee
which is contradiction, and the proof is finished in second case.

Next, we prove the sharpness of $p_N(u)$. It suffices to construct a sequence $\{u_n\}\subset W^{1, N}_0(\Omega)$ and a function  $u\in W^{1, N}_0(\Omega)$ such that
   $$\int_\Omega F^{N}(\nabla u_n)dx=1,\quad  u_n\rightharpoonup u\not \equiv 0  \quad \mbox{in} \quad  W^{1, N}_0(\Omega),$$
   $$\bigg(\int_\Omega F^{N}(\nabla u)dx\bigg)^{\frac{1}{N}}=\delta<1 \quad \mbox{and}\quad \int_{\Omega}\frac{e^{\lambda_{N}(1-\frac{\beta}{N}) (1-\delta^N)^{-\frac{1}{N-1}}|u_n|^{\frac{N}{N-1}}}}{F^{o}(x)^{\beta}}dx\rightarrow +\infty.$$
For $n\in \mathbb{N}$, let $r>0$, we define
\begin{equation*}
  \omega_{n}(x)=
  \left\{  \begin{array}{l}
          \frac{1}{N}\kappa_{N}^{-\frac{1}{N}}n^\frac{N-1}{N},      ~~~~~~0\leq F^{o}(x) \leq re^{-\frac{n}{N}},\\
          \kappa_{N}^{-\frac{1}{N}}\log(r/F^{o}(x))n^{-\frac{1}{N}},  ~~~~~re^{-\frac{n}{N}}\leq F^{o}(x) \leq r, \\
         0,  ~~~~F^{o}(x) \geq r.
         \end{array}
   \right.
 \end{equation*}
A straightforward calculation yields
 $$\int_{\Omega}F^{N}(\nabla \omega_n)dx=1,\quad  w_n\rightharpoonup 0 \quad  \mbox{in} \quad  W^{1, N}_0(\Omega).$$
Set $R=3r$, define
 \begin{equation*}
  u(x)=
  \left\{  \begin{array}{l}
          A,      ~~~~~~0\leq F^{o}(x) \leq \frac{2}{3}R,\\
          3A-\frac{3A}{R}F^{o}(x),  ~~~~~\frac{2}{3}R\leq F^{o}(x) \leq R, \\
         0,  ~~~~F^{o}(x) \geq R,
         \end{array}
   \right.
 \end{equation*}
where $A$ is a positive constant to be chosen in such a way that $(\int_{\Omega}F^{N}(\nabla u)dx)^{\frac{1}{N}}=\delta<1$. Denote $\mathcal{W}(R)=\{x\in \mathbb{R}^{N} : F^{o}(x)\leq R \}$ be a Wulff ball centered at the origin. Let $u_n=u+(1-\delta^N)^{\frac{1}{N}}\omega_{n}$, since $\nabla u$ and $\nabla \omega_n$ have disjoint supports, we have
$$\int_\Omega F^{N}(\nabla u_n)dx=\int_{\mathcal{W}(R)} F^{N}(\nabla u)dx+(1-\delta^N)\int_{\mathcal{W}(R)} F^{N}(\nabla w_n)dx=1$$
and $u_n\rightharpoonup u$ in  $W^{1, N}_0(\Omega)$. Thus
\begin{equation*}
\left.
\begin{aligned}[b]
             &\int_{\Omega}\frac{e^{\lambda_{N}(1-\frac{\beta}{N})(1-\delta^N)^{-\frac{1}{N-1}} |u_n|^{\frac{N}{N-1}}}}{F^{o}(x)^{\beta}}dx\\
              \geq &   \int_{\mathcal{W}(re^{-\frac{n}{N}})}\frac{e^{\lambda_{N}(1-\frac{\beta}{N})(1-\delta^N)^{-\frac{1}{N-1}} |A+(1-\delta^N)^{\frac{1}{N}}\omega_{n}|^{\frac{N}{N-1}}}}{F^{o}(x)^{\beta}}dx\\
        =&\int_{\mathcal{W}(re^{-\frac{n}{N}})}\frac{e^{\lambda_{N}(1-\frac{\beta}{N}) |C+\omega_{n}|^{\frac{N}{N-1}}}}{F^{o}(x)^{\beta}}dx\\
        =&e^{{N^{\frac{N}{N-1}} \kappa_{N}^{\frac{1}{N-1}}(1-\frac{\beta}{N}) }(C+\frac{1}{N}\kappa_N^{-\frac{1}{N}}n^\frac{N-1}{N})^{\frac{N}{N-1}}}\int_{\mathcal{W}(re^{-\frac{n}{N}})}\frac{1}{F^{o}(x)^{\beta}}dx\\
          =&e^{{N^{\frac{N}{N-1}} \kappa_{N}^{\frac{1}{N-1}}(1-\frac{\beta}{N}) }(C+\frac{1}{N}\kappa_N^{-\frac{1}{N}}n^\frac{N-1}{N})^{\frac{N}{N-1}}}\cdot \frac{N\kappa_N}{N-\beta}\rho^{N-\beta}|_0^{re^{-\frac{n}{N}}} \\
        \geq &e^{[C_1+((1-\frac{\beta}{N})n)^{\frac{N-1}{N}}]^{\frac{N}{N-1}}} \cdot \frac{N\kappa_N}{N-\beta}r^{N-\beta}e^{-(1-\frac{\beta}{N})n}  \\
          \geq &C_2 e^{[C_1+((1-\frac{\beta}{N})n)^{\frac{N-1}{N}}]^{\frac{N}{N-1}}} e^{-(1-\frac{\beta}{N})n}\rightarrow +\infty(n\rightarrow +\infty)
\end{aligned}
\right.
\end{equation*}
where $0\leq \beta <N$ and $C, C_1 , C_2$   are positive constants.   $\hfill\Box$

 \section{Lions type concentration-compactness principle in $\mathbb{R}^{N}$}
  As the similar procedure in Theorem 1.2, we can immediately get Theorem 1.3.  \\
 \textbf{Proof of Theorem 1.3.} Since
   $$\int_{\mathbb{R}^{N}}(F^{N}(\nabla u)+|u|^N)dx \leq \liminf_{n\rightarrow \infty}\int_{\mathbb{R}^{N}}(F^{N}(\nabla u_n)+|u_n|^N)dx =1.$$
We discuss it in two cases.

  Case 1: Let $0<\int_{\mathbb{R}^{N}}(F^{N}(\nabla u)+|u|^N)dx <1 $, we give the proof by contradiction. Assume that there exists some $p_1 < \bar{p}_N(u)$  and a subsequence of $\{u_n\}$(still denote $u_n$) such that
     \be \sup_n\int_{\mathbb{R}^{N}}\frac{\Phi(\lambda_{N}(1-\frac{\beta}{N})p_1 |u_n|^{\frac{N}{N-1}})}{F^{o}(x)^{\beta}} dx=+\infty.  \ee
Set $\Omega_L^{n}=\{x\in \mathbb{R}^{N} : |u_n(x)|\geq L\}$, where $L$ is a positive constant.
Since $0\leq \beta <N$, we have
 \begin{equation}
\left.
\begin{aligned}[b]
              &\int_{\mathbb{R}^{N}}\frac{\Phi(\lambda_{N}(1-\frac{\beta}{N})p_1 |u_n|^{\frac{N}{N-1}})}{F^{o}(x)^{\beta}}dx\\
              =& \int_{\Omega_L^{n}}\frac{\Phi(\lambda_{N}(1-\frac{\beta}{N})p_1 |u_n|^{\frac{N}{N-1}})}{F^{o}(x)^{\beta}}dx+\int_{\mathbb{R}^{N}\backslash\Omega_L^{n}}\frac{\Phi(\lambda_{N}(1-\frac{\beta}{N})p_1 |u_n|^{\frac{N}{N-1}})}{F^{o}(x)^{\beta}}dx\\
             \leq  & \int_{\Omega_L^{n}}\frac{\Phi(\lambda_{N}(1-\frac{\beta}{N})p_1 |u_n|^{\frac{N}{N-1}})}{F^{o}(x)^{\beta}}dx+C\int_{\mathbb{R}^{N}\backslash\Omega_L^{n}}\frac{|u_n|^N}{F^{o}(x)^{\beta}}dx\\
             \leq  & \int_{\Omega_L^{n}}\frac{\Phi(\lambda_{N}(1-\frac{\beta}{N})p_1 |u_n|^{\frac{N}{N-1}})}{F^{o}(x)^{\beta}}dx+C\int_{
             F^o(x)\leq 1}\frac{1}{F^{o}(x)^{\beta}}dx+C\int_{
             F^o(x)>1}|u_n|^Ndx\\
              \leq  & \int_{\Omega_L^{n}}\frac{\Phi(\lambda_{N}(1-\frac{\beta}{N})p_1 |u_n|^{\frac{N}{N-1}})}{F^{o}(x)^{\beta}}dx+C(p_1, L, N, \beta),
\end{aligned}
\right.
\end{equation}
and then
\be \sup_n\int_{\Omega_L^{n}}\frac{\Phi(\lambda_{N}(1-\frac{\beta}{N})p_1 |u_n|^{\frac{N}{N-1}})}{F^{o}(x)^{\beta}}dx=+\infty.  \ee
Let $v_n=u_n-L$, for any $\epsilon > 0$, we have
     \be  |u_n|^{\frac{N}{N-1}}\leq (1+\epsilon)v_n^{\frac{N}{N-1}}+C_\epsilon L^{\frac{N}{N-1}}.  \ee
Notice that
\begin{equation}
\left.
\begin{aligned}[b]
              \int_{\Omega_L^{n}}\frac{\Phi(\lambda_{N}(1-\frac{\beta}{N})p_1 |u_n|^{\frac{N}{N-1}})}{F^{o}(x)^{\beta}}dx
             \leq  & e^{\lambda_{N}(1-\frac{\beta}{N})p_1 C _\epsilon L^{\frac{N}{N-1}}}\int_{\Omega_L^{n}}\frac{e^{(1+\epsilon)\lambda_{N}(1-\frac{\beta}{N})p_1 |v_n|^{\frac{N}{N-1}}}}{F^{o}(x)^{\beta}}dx.
\end{aligned}
\right.
\end{equation}
Thus
\be \sup_n\int_{\Omega_L^{n}}\frac{e^{\lambda_{N}(1-\frac{\beta}{N})\overline{p}_1 |v_n|^{\frac{N}{N-1}}}}{F^{o}(x)^{\beta}}dx=\sup_n\int_{\Omega_L^{n}}\frac{e^{\lambda_{N}(1-\frac{\beta}{N})\overline{p}_1 |v_n|^{\frac{N}{N-1}}}}{F^{o}(x)^{\beta}}dx=+\infty,  \lb{4.1}\ee
where $\overline{p}_1=(1+\epsilon)p_1<\bar{p}_N(u)$. Now, we define
  $$T^L(u)=\min \{L, |u|\}sign(u) ~~~ and ~~~T_L(u)=u-T^L(u)$$
and choose $L$ so large that
     \be \frac{1-\int_{\mathbb{R}^{N}} (F^{N}(\nabla u)+|u|^N)dx}{1-\int_{\mathbb{R}^{N}} (F^{N}(\nabla T^L(u))+| T^L(u)|^N)dx} > \bigg(\frac{\overline{p}_1}{\bar{p}_N(u)}\bigg)^{N-1}. \lb{4.7}\ee
Since $T^L(u_n)$ is  bounded in $W^{1, N}({\mathbb{R}^{N}})$, hence, up to a subsequence, $T^L(u_n)\rightharpoonup T^L(u)$ in $W^{1, N}({\mathbb{R}^{N}})$ and
$T^L(u_n)\rightarrow T^L(u)$ a. e. in $\mathbb{R}^{N}$. Combing \eqref{4.1} and \eqref{(1.5)}, up to a subsequence, we have
 $$ \limsup_{n\rightarrow \infty}\int_{\Omega_L^{n}}F^{N}(\overline{p}_1^{\frac{N-1}{N}}\nabla v_n)dx\geq 1.  $$
Thus
 \be \int_{\Omega_L^{n}}F^{N}(\nabla v_n)dx= \int_{\Omega}F^{N}(\nabla T_L(u_n))dx\geq \bigg(\frac{1}{\overline{p}_1}\bigg)^{N-1}+o_n(1).   \ee
Then we have
\begin{equation*}
\left.
\begin{aligned}[b]
             &\bigg(\frac{1}{\overline{p}_1}\bigg)^{N-1}+\int_{\mathbb{R}^{N}}F^{N}(\nabla T^L(u_n))dx+\int_{\mathbb{R}^{N}}| T^L(u_n)|^Ndx+o_n(1)\\
              \leq&\bigg(\frac{1}{\overline{p}_1}\bigg)^{N-1}+\int_{\mathbb{R}^{N}}F^{N}(\nabla T^L(u_n))dx+\int_{\mathbb{R}^{N}}| u_n|^Ndx+o_n(1)\\
             \leq  & \int_{\mathbb{R}^{N}}F^{N}(\nabla T_L(u_n))dx+\int_{\mathbb{R}^{N} \backslash \Omega_L^{n}}F^{N}(\nabla u_n)dx+\int_{\mathbb{R}^{N}}| u_n|^Ndx\\
             =  & \int_{\Omega_L^{n}}F^{N}(\nabla u_n)dx+\int_{\mathbb{R}^{N}\backslash \Omega_L^{n}}F^{N}(\nabla u_n)dx+\int_{\mathbb{R}^{N}}| u_n|^Ndx=1.
\end{aligned}
\right.
\end{equation*}
From \eqref{4.7}, it holds
\begin{equation*}
\left.
\begin{aligned}[b]
             \overline{p}_1\geq & \frac{1}{(1- \liminf_{n\rightarrow \infty}\int_{\mathbb{R}^{N}} (F^{N}(\nabla T^L(u_n))+| T^L(u_n)|^N)dx )^{\frac{1}{N-1}}}\\
            \geq & \frac{1}{(1- \int_{\mathbb{R}^{N}} (F^{N}(\nabla T^L(u))+| T^L(u)|^N)dx )^{\frac{1}{N-1}}}\\
             > & \frac{\overline{p}_1}{\bar{p}_N(u)}   \frac{1}{(1- \int_{\mathbb{R}^{N}}(F^{N}(\nabla u)+|u|^N)dx )^{\frac{1}{N-1}}}=\overline{p}_1,
\end{aligned}
\right.
\end{equation*}
which is a contradiction.  The proof is finished in the first case.

 Case 2: Let $\int_{\mathbb{R}^{N}}(F^N (\nabla u)+|u|^N)dx=1 $.  Since $u_n \rightharpoonup u$  and $W^{1, N}(\mathbb{R}^{N})$ is a uniformly convex
Banach space, we know that $u_n\rightarrow u$ in $W^{1, N}(\mathbb{R}^{N})$. Thus, by Proposition 1 in [16], there exists some $v\in W^{1, N}(\mathbb{R}^{N})$, such that up to
a subsequence, $|u_n(x)|\leq v(x)$ a.e. in $W^{1, N}(\mathbb{R}^{N})$. Denote
 $$D=\{x\in\mathbb{R}^{N}:\int_{\mathbb{R}^{N}} F^N(\nabla v)dx \leq 1, v(x)>1\},$$
 we have
\begin{equation*}
\left.
\begin{aligned}[b]
  \int_{\mathbb{R}^{N}\backslash D}\frac{\Phi(\lambda_{N}(1-\frac{\beta}{N})p_1 |v|^{\frac{N}{N-1}})}{F^{o}(x)^{\beta}} dx\leq C(p_1, N, \beta).
\end{aligned}
\right.
\end{equation*}
Indeed,
\begin{equation*}
\left.
\begin{aligned}[b]
  &\int_{\mathbb{R}^{N}\backslash D}\frac{\Phi(\lambda_{N}(1-\frac{\beta}{N})p_1 |v|^{\frac{N}{N-1}})}{F^{o}(x)^{\beta}} dx\\
  &\leq  \int_{\{v(x)>1\}}\frac{1}{F^{o}(x)^{\beta}} \sum_{k=N-1}^{\infty} \frac{[\lambda_{N}(1-\frac{\beta}{N}p_1]^k|u|^{\frac{kN}{N-1
  }}}{k !} dx\\
   &\leq  \int_{\{v(x)> 1\}}\frac{1}{F^{o}(x)^{\beta}} \sum_{k=N-1}^{\infty} \frac{[\lambda_N(1-\frac{\beta}{N})p_1]^k|u|^{N}}{k !} dx\\
   &\leq  \int_{\{F^o(x)\geq 1\}} \sum_{k=N-1}^{\infty} \frac{[\lambda_N(1-\frac{\beta}{N})p_1]^k|u|^{N}}{k !} dx+  \int_{\{F^o(x)<1\}}\frac{1}{F^{o}(x)^{\beta}}  \sum_{k=N-1}^{\infty} \frac{[\lambda_N(1-\frac{\beta}{N})p_1]^k}{k !} dx\\
   &\leq C(p_1, N, \beta)
\end{aligned}
\right.
\end{equation*}
From Lemma 2.4 and \eqref{2.5}, we have
\begin{equation}
\left.
\begin{aligned}[b]
&\int_{\mathbb{R}^{N}}\frac{\Phi(\lambda_{N}(1-\frac{\beta}{N})p_1 |u_n|^{\frac{N}{N-1}})}{F^{o}(x)^{\beta}}dx\\
\leq&\int_{\mathbb{R}^{N}}\frac{\Phi(\lambda_{N}(1-\frac{\beta}{N})p_1 |v|^{\frac{N}{N-1}})}{F^{o}(x)^{\beta}} dx\\
  \leq&\int_{\mathbb{R}^{N}\backslash D}\frac{\Phi(\lambda_{N}(1-\frac{\beta}{N})p_1 |v|^{\frac{N}{N-1}})}{F^{o}(x)^{\beta}} dx+\int_{ D}\frac{\Phi(\lambda_{N}(1-\frac{\beta}{N})p_1 |v|^{\frac{N}{N-1}})}{F^{o}(x)^{\beta}} dx\\
  \leq&C(p_1, N, \beta)+\int_{\mathcal{W}(R)}\bigg(\frac{1}{(F^o)^\beta}\bigg)^\star(x)\cdot\Phi(\lambda_{N}(1-\frac{\beta}{N})p_1 |v^{\star}|^{\frac{N}{N-1}}) dx\\
   \leq&C(p_1, N, \beta)+\int_{\mathcal{W}(R)}\frac{e^{\lambda_{N}(1-\frac{\beta}{N})p_1 |v^{\star}|^{\frac{N}{N-1}}}}{F^{o}(x)^{\beta}}  dx
\end{aligned}
\right.  \lb{4.9}
\end{equation}
where $\mathcal{W}(R)=\{x\in \mathbb{R}^{N} : F^{o}(x)\leq R \}$ be a Wulff ball and $|\mathcal{W}(R)|=|D|$. We know that $\int_{\mathcal{W}(R)}F^N(\nabla u^{\star})dx\leq\int_{\mathcal{W}(R)}F^N (\nabla u)dx$ by Lemma 2.3. Hence, the result follows from \eqref{(1.5)}.

Next, we prove the sharpness of $\bar{p}_N(u)$. It suffices to construct a sequence $\{u_n\}\subset W^{1, N}(\mathbb{R}^{N})$ and a function  $u\in W^{1, N}(\mathbb{R}^{N})$ such that
   $$\int_{\mathbb{R}^{N}}( F^{N}(\nabla u_n)+|u_n|^N)dx=1,\quad  u_n\rightharpoonup u\not \equiv 0  \quad \mbox{in} \quad  W^{1, N}(\mathbb{R}^{N}),$$
   $$\bigg(\int_{\mathbb{R}^{N}}( F^{N}(\nabla u_n)+|u_n|^N)dx\bigg)^{\frac{1}{N}}=\delta<1$$ and
   $$ \int_{\mathbb{R}^{N}}\frac{\Phi({\lambda_{N}(1-\frac{\beta}{N}) \bar{p}_N(u)|u_n|^{\frac{N}{N-1}}})}{F^{o}(x)^{\beta}}dx\rightarrow +\infty.$$
For $n\in \mathbb{N}$, let $r>0$, we define
\begin{equation*}
  \omega_{n}(x)=
  \left\{  \begin{array}{l}
          \frac{1}{N}\kappa_{N}^{-\frac{1}{N}}n^\frac{N-1}{N},      ~~~~~~0\leq F^{o}(x) \leq re^{-\frac{n}{N}},\\
          \kappa_{N}^{-\frac{1}{N}}\log(r/F^{o}(x))n^{-\frac{1}{N}},  ~~~~~re^{-\frac{n}{N}}\leq F^{o}(x) \leq r, \\
         0,  ~~~~F^{o}(x) \geq r.
         \end{array}
   \right.
 \end{equation*}
A straightforward calculation yields
 $$ w_n\rightharpoonup 0 \quad  \mbox{in} \quad  W^{1, N}_0(\Omega),  \quad \int_{\mathbb{R}^{N}}F^{N}(\nabla \omega_n)dx=1,\quad \int_{\mathbb{R}^{N}} |\omega_n|^Ndx\rightarrow 0.$$
Set $R=3r$, define
 \begin{equation*}
  u(x)=
  \left\{  \begin{array}{l}
          A,      ~~~~~~0\leq F^{o}(x) \leq \frac{2}{3}R,\\
          3A-\frac{3A}{R}F^{o}(x),  ~~~~~\frac{2}{3}R\leq F^{o}(x) \leq R, \\
         0,  ~~~~F^{o}(x) \geq R,
         \end{array}
   \right.
 \end{equation*}
where $A$ is a positive constant to be chosen in such a way that $(\int_{\mathbb{R}^{N}}( F^{N}(\nabla u)+|u|^N)dx)^{\frac{1}{N}}=\delta<1$. Denote $\mathcal{W}(R)=\{x\in \mathbb{R}^{N} : F^{o}(x)\leq R \}$ be a Wulff ball centered at the origin. Set $v_n=u+(1-\delta^N)^{\frac{1}{N}}\omega_{n}$, we have
\begin{equation*}
\left.
\begin{aligned}[b]
&\int_{\mathbb{R}^{N}} F^{N}(\nabla v_n)dx\\
=&\int_{\mathcal{W}(R)} F^{N}(\nabla u)dx+(1-\delta^N)\int_{\mathcal{W}(R)} F^{N}(\nabla w_n)dx\\
=&\int_{\mathcal{W}(R)} F^{N}(\nabla u)dx+(1-\delta^N).
\end{aligned}
\right.
\end{equation*}
Moreover, we have
\begin{equation*}
\left.
\begin{aligned}[b]
             \int_{\mathbb{R}^{N}} |v_n|^Ndx=&\int_{\mathbb{R}^{N}} |u+(1-\delta^N)^{\frac{1}{N}}\omega_{n}|^Ndx\\
        =&\int_{\mathbb{R}^{N}} |u|^Ndx+r_n,
\end{aligned}
\right.
\end{equation*}
where $r_n=O(n^{-\frac{1}{N}})$ as $n\rightarrow +\infty$. Thus we have $\int_{\mathbb{R}^{N}}( F^{N}(\nabla v)+|v|^N)dx=1+r_n$. Let $u_n=\frac{v_n}{(1+r_n)^{\frac{1}{N}}}$,
it holds
$$ \int_{\mathbb{R}^{N}}( F^{N}(\nabla u_n)+|u_n|^N)dx=1  \quad  u_n\rightharpoonup u \quad  \mbox{in} \quad  W^{1, N}(\mathbb{R}^{N}).$$
Then
\begin{equation*}
\left.
\begin{aligned}[b]
             &\int_{\mathbb{R}^{N}}\frac{\Phi({\lambda_{N}(1-\frac{\beta}{N}) \bar{p}_N(u)|u_n|^{\frac{N}{N-1}}})}{F^{o}(x)^{\beta}}dx\\
             \geq &   \int_{\mathcal{W}(re^{-\frac{n}{N}})}\frac{e^{\lambda_{N}(1-\frac{\beta}{N})(1-\delta^N)^{-\frac{1}{N-1}} |u_n|^{\frac{N}{N-1}}}}{F^{o}(x)^{\beta}}dx+C(u)\\
                         \geq &   \int_{\mathcal{W}(re^{-\frac{n}{N}})}\frac{e^{\lambda_{N}(1-\frac{\beta}{N})((1+r_n)(1-\delta^N))^{-\frac{1}{N-1}} |A+(1-\delta^N)^{\frac{1}{N}}\omega_{n}|^{\frac{N}{N-1}}}}{F^{o}(x)^{\beta}}dx+C(u)\\
        =&\int_{\mathcal{W}(re^{-\frac{n}{N}})}\frac{e^{\lambda_{N}(1-\frac{\beta}{N})(1+r_n)^{-\frac{1}{N-1}}  |C+\omega_{n}|^{\frac{N}{N-1}}}}{F^{o}(x)^{\beta}}dx+C(u)\\
        =&e^{{N^{\frac{N}{N-1}} \kappa_{N}^{\frac{1}{N-1}}(1-\frac{\beta}{N})(1+r_n)^{-\frac{1}{N-1}}  }(C+\frac{1}{N}\kappa_N^{-\frac{1}{N}}n^\frac{N-1}{N})^{\frac{N}{N-1}}}\int_{\mathcal{W}(re^{-\frac{n}{N}})}\frac{1}{F^{o}(x)^{\beta}}dx+C(u)\\
          =&e^{{N^{\frac{N}{N-1}} \kappa_{N}^{\frac{1}{N-1}}(1-\frac{\beta}{N}) }(1+r_n)^{-\frac{1}{N-1}} (C+\frac{1}{N}\kappa_N^{-\frac{1}{N}}n^\frac{N-1}{N})^{\frac{N}{N-1}}}\cdot \frac{N\kappa_N}{N-\beta}\rho^{N-\beta}|_0^{re^{-\frac{n}{N}}} +C(u)\\
        \geq &e^{[C_1+(1+r_n)^{-\frac{1}{N}} ((1-\frac{\beta}{N})n)^{\frac{N-1}{N}}]^{\frac{N}{N-1}}} \cdot \frac{N\kappa_N}{N-\beta}r^{N-\beta}e^{-(1-\frac{\beta}{N})n} +C(u)\\
          \geq &C_2 e^{[C_1+(1+r_n)^{-\frac{1}{N}} ((1-\frac{\beta}{N})n)^{\frac{N-1}{N}}]^{\frac{N}{N-1}}} e^{-(1-\frac{\beta}{N})n}+C(u)\rightarrow +\infty(n\rightarrow +\infty)
\end{aligned}
\right.
\end{equation*}
where $0\leq \beta <N$ and $C, C_1 , C_2$   are positive constants.   $\hfill\Box$

\def\refname{References }

\end{document}